# Application of Computational Tools for Finitely Presented Groups

## GEORGE HAVAS AND EDMUND F. ROBERTSON

ABSTRACT. Computer based techniques for recognizing finitely presented groups are quite powerful. Tools available for this purpose are outlined. They are available both in stand-alone programs and in more comprehensive systems.

A general computational approach for investigating finitely presented groups by way of quotients and subgroups is described and examples are presented. The techniques can provide detailed information about group structure. Under suitable circumstances a finitely presented group can be shown to be soluble and its complete derived series can be determined, using what is in effect a soluble quotient algorithm.

## 1. Introduction

We address the problem of determining structural information about a finitely presented group. Finitely presented groups arise naturally in many ways. Some of these are purely combinatorial while others originate topologically. When such groups occur their structures may not be known. We look at some reasonably standard approaches and show how effective they can be with readily available computer programs on workstation size machines.

We give an overview of tools available for computing with finitely presented groups. We include both stand-alone implementations and also packages designed for computational group theory. Some further information on algorithms for groups is included in [**14**].

A good introduction to group theory is provided by Rotman [**43**]. Important concepts we use in describing the computations in this paper are the following.

1991 *Mathematics Subject Classification.* Primary 20-04, 20F05; Secondary 20F14.

We are grateful for EC grant SC1*0003-C(EDB) which made the first author's visit to St Andrews possible and thus led to the writing of this paper. The computer facilities used are partially supported by SERC grant GR/E91127. The research of the first author was partially supported by ARC grant A49030651.

This paper is in final form and no version of it will be submitted for publication elsewhere.

A finite presentation $\langle X \mid R \rangle$ consists of a finite set $X$ of generators and a finite set $R$ of relators, which are words in the generators and their inverses. A group $G$ is finitely presented if there is a finite presentation $\langle X \mid R \rangle$ with $G \cong F/N$ where $F$ is the free group on $X$ and $N$ is the normal closure of $R$ in $F$. Given a group $G$ with subgroups $H$ and $K$ define $[H, K]$ to be the subgroup generated by all commutators $[h, k] = h^{-1}k^{-1}hk$ with $h \in H$ and $k \in K$. Then $G' = [G, G]$ is the derived group of $G$, and $G/G'$ is the largest abelian factor group of $G$. The derived series of $G$ is defined inductively by $G^{(0)} = G$ and $G^{(i+1)} = G^{(i)'}$. A soluble group $G$ is one for which $G^{(n)} = 1$ for some $n$ and the derived length of $G$ is the least such $n$. The lower central series of $G$ is defined inductively by $\gamma_1(G) = G$ and $\gamma_{i+1}(G) = [\gamma_i(G), G]$, and a nilpotent group is one for which $\gamma_n(G) = 1$ for some $n$.

Given a finitely presented group we can ask a number of questions about it. These include: is it trivial? is it finite? is it abelian? is it nilpotent? is it soluble? what is its structure?

These questions about finitely presented groups are unsolvable in general. Algorithms exist for verifying the truth of positive answers to the first four, but the last two are harder. Thus, subgroups of finitely presented soluble groups are not necessarily finitely generated. Furthermore, the structure question is more complicated in that it requires more than a yes/no answer, and different answers might well suffice for different purposes.

The verification algorithms lead to procedures which attempt to answer the questions, but which fail to terminate if the answer is negative. For example, given a finitely presented group which is finite, coset enumeration will (given adequate resources) produce a multiplication table for the group, from which the answers can be deduced. Alternatively the Knuth-Bendix procedure for strings may be used. Comprehensive discussions of these and related issues appear in the book by Sims [**46**].

Given a group presentation $G = \langle g_1, \dots, g_d \mid R_1, \dots, R_n \rangle$ the first thing we do to investigate the structure question is compute $G/G'$. If $G/G'$ is infinite then perhaps we already have enough information. We consider situations where $G/G'$ is finite.

Two apparently different cases arise: if $|G/G'| = 1$ either $G$ is trivial or we have an insoluble group; otherwise we have a soluble group (perhaps acting on something underneath). At this stage we can simply try one of the most important computational tools for finitely presented groups, coset enumeration. We can try a coset enumeration over the trivial subgroup. If the enumeration succeeds we know the group order and have a faithful permutation representation, so that we can readily enough answer questions about the group.

Let us assume that coset enumerations over both trivial and obvious cyclic subgroups fail. What next?

We proceed by looking at (nontrivial) quotients and subgroups of the given group. If $|G/G'| \neq 1$ then study of the derived series of $G$ is an obvious route.

The structure of this series may well tell us as much as we want to know about the group. On the other hand if $|G/G'| = 1$, we show how interesting quotients may be found and how we can investigate corresponding subgroups.

Our examples include cases where complete group structure is effectively elucidated, both for finite and infinite groups. All computations were done using stand-alone computer programs running on Sun workstations at St Andrews. At the time of writing (end of 1991) we found they could be done most conveniently with such stand-alone programs. We have updated this paper at final submission time (end of 1993) to include more recent references.

## 2. Computational tools

We indicate computational tools that may be used for answering questions about finitely presented groups. First we outline the tools used in our examples, then we mention other tools that are available.

**2.1. Coset enumeration.** Coset enumeration is the basis of important techniques for investigating finitely presented groups. Coset enumeration programs implement systematic procedures for enumerating the cosets of a subgroup $H$ of finite index in a group $G$, given defining relations for $G$ and words generating $H$. Computer implementations are based on methods initially described by Todd and Coxeter [**47**]. Some recent advances are described by Havas [**20**] and further improvements are available in implementations but otherwise unpublished. Variants of these implementations are used for the computations exemplified in this paper. Neubüser [**35**] gives many applications of coset enumeration. Sims [**46**], in a substantial chapter, gives a formal account of coset enumeration in terms of automata and proves interesting new results on coset enumeration behaviour. Linton [**33**] has implemented double coset enumeration and shown how it can be effectively applied in suitable situations. In other work [**34**] he shows how to use a version of coset enumeration to construct matrix representations for finitely presented groups.

**2.2. Subgroup presentation.** This topic is addressed by Neubüser [**35**]. Based on theorems of Reidemeister and Schreier, various programs for computing presentations of subgroups of finite index in finitely presented groups have been developed. These include those by Havas [**19**] and Arrell and Robertson [**1**], which are used in our examples.

**2.3. Tietze transformation.** Subgroup presentations which are produced using Reidemeister-Schreier processes generally have large numbers of generators and relators and are often poorly suited for human or computer work. A theorem of Teitze implies that, given two presentations for a group, a sequence of simple transformations exists which converts one presentation to the other. Even though there is no general algorithm for finding such a sequence, Tietze transformation programs which take a "bad" presentation and produce a "good" presentation have been developed by Havas, Kenne, Richardson and Robertson [**22**] and Robertson [**42**], and are used in the examples presented here.

**2.4. Abelian group recognition.** The theory of finitely generated abelian groups provides a method for completely classifying such groups. For finitely presented abelian groups a canonical form can be obtained by using algorithms for Smith normal form computation for integer matrices. When the presentations are "bad" (many generators and many relations) sophisticated methods may be required to successfully perform such calculations. Programs based on Havas and Sterling [**25**] are used for this purpose in our examples. More recent developments appear in [**21**].

**2.5. Nilpotent quotient calculation.** An algorithm which computes $p$-quotients of finitely presented groups is described by Havas and Newman [**23**]. A new implementation, developed by Newman and O'Brien at the Australian National University, Canberra, incorporates a better collector and uses automorphism methods to handle groups satisfying identical relations more efficiently. Presentations for certain extensions of $p$-groups can be constructed using the technique of $p$-group generation. This algorithm and its implementation are described in O'Brien [**40**]. Some more recent developments on $p$-quotient algorithm implementation are described by Celler, Newman, Nickel and Niemeyer in an ANU Research Report. A general nilpotent quotient program (with no dependency on a prime $p$) has been developed by Nickel and shows promise.

**2.6. Soluble quotient calculation.** A more difficult area is the computation of soluble quotients. Two recent PhD dissertations address this area with some success. Wegner [**48**] and Niemeyer [**39**] both implement algorithms and apply them to various examples.

**2.7. Low index subgroups.** Given a finitely presented group it is possible to search combinatorially for subgroups of low index by checking coset tables. Descriptions of such methods are given by Neubüser [**35**] and Sims [**46**].

**2.8. Quotients and kernels.** Holt and Rees [**29**] have developed an effective program **quotpic** for studying quotients of finitely presented groups and their kernels. They provide "at least potential access to all quotients up to order 10000". This involves many of the tools listed above. The program has an effective X-windows interface which gives a clear graphical representation of relevant subgroup lattices.

**2.9. Knuth-Bendix.** The Knuth-Bendix procedure [**30**] has gradually been playing a more important role in computational group theory. Early applications have been mainly with infinite groups ([**18**, **31**]). More recently Sims [**45**] has shown some situations where Knuth-Bendix outperforms coset enumeration. Epstein, Holt and Rees [**16**] have made significant progress on automatic groups.

**2.10. Isomorphism/nonisomorphism testing.** In a very neat combination of many of the above techniques Holt and Rees [**28**] have developed a program which endeavours to determine isomorphism or nonisomorphism of two finitely presented groups. Isomorphism testing relies on the Knuth-Bendix procedure for finding an explicit mapping while nonisomorphism testing relies on finding some structural difference. The program switches to and from each kind of testing hoping to discover isomorphism or nonisomorphism. In addition, O'Brien [**41**] has developed an algorithm for $p$-groups. Even though this problem is unsolvable in general, these methods show an attractive success rate.

**2.11. Implementations.** Major software development projects have led to three packages which include comprehensive facilities for computational group theory. Cayley [**8**], GAP [**44**] and MAGMA [**12**] include algorithms which implement most of the above tools, with higher level interfaces. Their substantial reference manuals indicate their capabilities. However, at the time of initial writing none of these packages could easily handle the type of calculations presented here in their entirety.

Subsets of these tools have been put together in various ways to give integrated computing facilities. Thus SPAS [**17**] provided access to coset enumeration, subgroup presentation, Tietze transformation, $p$-quotient and abelian group recognition programs via a simple command language. In St Andrews (see Heggie [**26**]) similar facilities are available loosely interconnected in the sense that the various stand-alone programs produce outputs which (possibly with the help of simple filters) are accepted as inputs for the others. This is the environment we used for our calculations. Indeed the **quotpic** package [**29**] now provides a very user-friendly way of addressing these kinds of problems.

**2.12. Applications.** There are too many applications of these kinds of tools to problems in finitely presented groups for us to attempt to list them all here. The conference proceedings Leech [**32**] and Atkinson [**2**] include earlier representative samples, while two special issues of the Journal of Symbolic Computation [**9, 11**] are more recent. Many of the Cayley applications listed by Cannon [**10**] use these methods.

## 3. The basic approach and its utilization

Given a finitely presented group whose structure is unknown we aim to obtain further information, perhaps complete structural information. The approach we take is not new. However we show that quite strong results may be easily achieved with readily available computing facilities. All of the examples presented here were resolved with stand-alone computer programs for coset enumeration, abelian group recognition, Tietze transformation and subgroup presentation, as available at St Andrews.

Consider the following situation. We have a group $G$ with unknown structure, given by a finite presentation $G = \langle g_1, \ldots, g_d \mid R_1, \ldots, R_n \rangle$. We compute $G/G'$ and discover that $G/G'$ is finite. This computation can often be readily done by hand. However more difficult cases arise, as we shall see, which are then resolved using an abelian group recognition program.

First consider the case that $G/G'$ is not trivial. Then $G$ has a nontrivial soluble quotient. One general approach is to try to follow the derived series of $G$ to find $G/G'$, $G'/G''$, … , or perhaps some other subnormal series. This kind of approach has been successfully followed elsewhere, including Neubüser and Sidki [**36**] and Newman and O'Brien [**38**].

When $G/G'$ is nontrivial it gives us not only the first step in the derived series but also the first step in the lower central series. Indeed applications arise where nilpotent quotients of $G$ tell us as much as we want to know, and $p$-quotients are readily computable. A striking example of the successful use of this type of approach by computer appears in the investigation of the Fibonacci group $F(2, 9)$ by Havas, Richardson and Sterling [**24**]; Newman [**37**], with a sophisticated application of the Golod-Šaferevič theorem, went on to prove that $F(2, 9)$ is infinite. We demonstrate the general method in Example 1, where we determine the complete derived series of an infinite group, showing that the approach yields an effective soluble quotient algorithm in suitable circumstances. Furthermore, in our particular case, we prove that the group is soluble.

When $G/G'$ is trivial, $G$ is perfect. The study of perfect groups is a major task (see Holt and Plesken[**27**]). There are no immediately obvious finite quotients to consider. Suffice it to say that finite quotients may exist which are (direct products of nonabelian) simple groups.

Various techniques exist for finding such quotients. There are low index subgroup programs which can sometimes find such quotients and have been used to do so. Alternatively, random coincidence procedures (such as those used by Cannon and Havas [**13**]) are also often effective. The method we use here is to find quotients by adding suitable relations to the given presentation.

A look at the two generator presentations for the simple groups of order less than one million (see Cannon, McKay and Young [**15**], Campbell and Robertson [**4, 5**]) shows that low powers of short words (like simple products or commutators of generators) are often trivial in the groups. So, to find a finite quotient we systematically add low powers of such words to given presentations and try coset enumerations. We often find nontrivial quotients this way. If successful, this shows that $G$ itself is not trivial, which is quite possibly new and interesting information in itself.

Having found a nontrivial quotient we go further. We look at the subgroup, $H$, which is the kernel of the homomorphism from $G$ to the quotient. Since we have a presentation for $G$ and a coset table for $G/H$ we can use a subgroup presentation (Reidemeister-Schreier) program to obtain a presentation for the subgroup $H$. Since the number of Schreier generators is $|G/H|(d - 1) + 1$ and the number of Reidemeister relators is $|G/H|n$, it is generally convenient (and sometimes necessary) for further computations that we simplify the ensuing presentations.

We do this with a Tietze transformation program.

Now look at where we are. If we do not recognize $H$ then we have a finite presentation for a group with unknown structure, now $H$. So we start again at the beginning. (This is one of the ways that finitely presented groups arise combinatorially. Sometimes the presentations arising this way are bad enough to require a computer to enable us to recognize even their abelian quotients.)

We give two examples where this approach is successfully followed to tell us all that we want to know about the groups involved.

Observe that the distinction that we have made so far between the two cases, $G/G'$ trivial and $G/G'$ nontrivial, is not one of principle. Basically what we do in both cases is find a nontrivial finite quotient $G/H$ of $G$ and look at the kernel of the homomorphism $G \to G/H$. When $G/G'$ is itself nontrivial, $G'$ is an immediate candidate for the kernel. Even when $G/G'$ is nontrivial we sometimes successfully find alternative subgroups to study, particularly when $|G/G'|$ is so large that direct computations with $G'$ are difficult. The key to the success of the method is the existence of useful subgroups of small enough index that we can discover and compute with. Our examples show that index well into the hundreds can be within range, sometimes directly and sometimes by repeating the process.

## 4. Three applications

**4.1. Example 1.** $G(2) = \langle a, b \mid a^6 = b^6 = 1, \ ab^2 = ba^2 \rangle$. The underlying family of groups $G(n)$ with generators of odd order is studied by Campbell and Robertson [**6**]. The family with generators of even order is considered by Campbell, Robertson and Thomas [**7**].

The computations:

   (i) Abelianize $G$ by adding $[a, b] = 1$.
  (ii) Find a coset table $T$ for $G/G' \cong C_6$.
 (iii) Find a presentation for $G'$ by running a subgroup presentation program on relators for $G$ and the coset table $T$.
  (iv) Apply a Tietze transformation program to simplify the presentation for $G'$ (to a presentation on 2 generators).
   (v) Repeat steps 1 to 3, with $G'$ instead of $G$. This yields $G'/G'' \cong C_4^2$, and a presentation for $G''$. The Tietze transformation program then shows exactly what $G''$ is. It (automatically) gives a final presentation $G'' = \langle r, s, t \mid [r, s] = [r, t] = [s, t] = 1 \rangle$. Thus $G'' \cong C_\infty^3$.

This shows that $G$ is soluble of derived length 3, with derived factors $C_6$, $C_4^2$ and $C_\infty^3$. We can readily redo the calculation with specific subgroup generators to get more details.

Let $G' = \langle x, y \rangle$ where $x = aba^{-1}b^{-1}$ and $y = ab^{-1}a^{-1}b$. We find the action of $a$ and $b$ on $x$ and $y$ to be $x^a = x^{-1}y^{-1}x$, $x^b = y^{-1}$, $y^a = x^{-1}y$, $y^b = y^{-1}x^{-1}y^2$. We do this by using the Reidemeister-Schreier method to rewrite $x$, $y$ and their conjugates in terms of the Schreier generators. Then we use the Tietze transformation program to eliminate the Schreier generators, presenting $G'$ on $x$ and $y$ and finding expressions for the conjugates (which appear as redundant subgroup generators) in terms of $x$ and $y$.

Similarly, let $G'' = \langle r, s, t \rangle$ where $r = xyx^{-1}y^{-1}$, $s = xy^{-1}x^{-1}y$, $t = x^{-1}y^{-1}xy$. Then we find $r^x = r^{-1}st^{-1}$, $r^y = s^{-1}$, $s^x = t^{-1}$, $s^y = r^{-1}$, $t^x = s^{-1}$, $t^y = r^{-1}st^{-1}$.

All of these calculations could in fact be eventually accomplished by hand. A similar computation with $G(3) = \langle a, b, c \mid a^6 = b^6 = c^6 = 1,\ ab^2 = ba^2,\ ac^2 = ca^2,\ bc^2 = cb^2 \rangle$ shows that $G(3)$ is also a soluble group of derived length 3, this time with derived factors $C_6, C_4^4$ and $C_\infty^9$. We suggest that this calculation is beyond the range of hand work, involving a subgroup of index 256 with nine generators.

**4.2. Example 2.** $G = \langle a, b \mid a^3 = b^5 = (aba^{-1}b^{-1}ab^2)^2 = 1 \rangle$. This is a generalized triangle group which "just fails" the infiniteness condition of Baumslag, Morgan and Shalen [**3**] for such groups and has been studied by Levin and Rosenberger. $G$ is perfect. It is easy to find finite homomorphic images. We systematically add relations of the form $(ab)^n = 1$ and try coset enumerations.

If we add $(ab)^2 = 1$ to the presentation for $G$, the resulting group is trivial, which is easily shown by coset enumeration over the trivial subgroup. Addition of the relation $(ab)^3 = 1$ gives $A_5$ as an image. The kernel of this map is perfect and is given to us as an unpleasantly presented 11 generator, 12 relator group, so we continue on. If we add $(ab)^4 = 1$, coset enumerations over the trivial subgroup do not complete with reasonable space restrictions, and adding extra relations leads to total collapse.

Now we add $(ab)^5 = 1$ to the presentation for $G$. Coset enumerations fail to complete, so we add another relation. We try $[a, b]^n = 1$, for various $n$. Addition of $[a, b]^2 = 1$ leads to total collapse. Addition of $[a, b]^3 = 1$ is useful. The image has order 1920, which is too large an index for us to proceed comfortably to subgroup presentation. To reduce the index, we take the subgroup $\langle b \rangle$ instead of the trivial subgroup in this image and construct a coset table with 384 cosets. The subgroup presentation program gives us a presentation on 210 generators and 225 relators for the corresponding subgroup $H$, say, of $G$

We prefer to simplify this presentation for $H$ using the Tietze transformation program before applying the abelian group recognition program to identify its abelian quotient. In part this is because we may sometimes (as in Example 1) be able to identify $H$ itself, rather than just its abelian quotient. In this case we quickly eliminate 120 generators without much increase in presentation length. It is clear by this stage that further Tietze eliminations will become increasingly difficult, so we look at the abelian quotient of $H$ as presented on these 90 generators. Determinant calculations show that the order of the torsion subgroup of the abelian quotient of $H$ is a divisor of $2^{13} \times 5^5 \times 3011$. Then modular diagonalization of the relation matrix gives invariants $C_2^4 \times C_4^2 \times C_8 \times C_\infty^{14}$ mod $2^{13}$, $C_5^3 \times C_\infty^{14}$ mod $5^5$ and $C_\infty^{14}$ mod 3011.

This is more than enough to tell us that $G$ is infinite.

**4.3. Example 3.** $G = \left\langle a, b \mid a^3 = b^{31} = (aba^{-1}b^2)^2 = (ab^2a^{-1}b)^2 = (ab^3a^{-1}b^{-11})^2 = (ab^4a^{-1}b^{13})^2 = 1 \right\rangle$. This presentation arises during the study of efficient presentations for direct products of simple groups.

$G/G' \cong C_3$ so we apply the method of Example 1 to find a presentation for $G'$. We get a 3 generator, 15 relator presentation. Let $x$, $y$, $z$ be the generators. A coset enumeration in $G'$ over $H = \left\langle x, y \right\rangle$ gives $|G' : H| = 32736$. There are six relations in the presentation for $H$ which involve only the generators $x$ and $y$, namely $x^{31} = y^{31} = (x^2y)^2 = (xy^2)^2 = (x^3y^{-11})^2 = (x^4y^{13})^2 = 1$. If $\bar{H}$ is the group presented by these two generators and six relations then $\bar{H}$ is a perfect group and, since $|\bar{H} : \left\langle x \right\rangle| = 1056$ (by coset enumeration), $|\bar{H}| = 32736$. Now, by Holt and Plesken [**27**], $\bar{H} \cong SL(2, 32)$ and a little more work shows that $G' \cong SL(2, 32) \times SL(2, 32)$ and that $G \cong C_3 \times SL(2, 32) \times SL(2, 32)$.

Key Centre for Software Technology, Department of Computer Science, University of Queensland, Queensland 4072, Australia
    *E-mail address*: havas@cs.uq.oz.au

Department of Mathematical and Computational Sciences, University of St Andrews, St Andrews KY16 9SS, Scotland
    *E-mail address*: efr@st-andrews.ac.uk